\documentclass[12pt,letterpaper]{article}

\usepackage{epsfig}
\usepackage{amsmath}
\usepackage{amssymb}
\usepackage{amsthm}
\usepackage{indentfirst}
\usepackage{xspace}
\usepackage{multirow}
\usepackage{hyperref}
\usepackage{verbatim}
\usepackage[letterpaper,margin=1in,headheight=15pt]{geometry}
\usepackage{mdwlist}
\usepackage{framed}
\usepackage{color}

\definecolor{shadecolor}{rgb}{0.85,0.85,0.85}
\definecolor{darkred}{rgb}{0.5,0.15,0.15}
\hypersetup{colorlinks=true,urlcolor=darkred,linkcolor=darkred,citecolor=darkred}

\numberwithin{equation}{section}

\newcommand{\cB}{\ensuremath{\mathcal B}}

\newcommand{\cD}{\ensuremath{\mathcal D}}
\newcommand{\cG}{\ensuremath{\mathcal G}}

\newcommand{\cK}{\ensuremath{\mathcal K}}
\newcommand{\cZ}{\ensuremath{\mathcal Z}}
\newcommand{\cM}{\ensuremath{\mathcal M}}
\newcommand{\cO}{\ensuremath{\mathcal O}}

\newcommand{\cX}{\ensuremath{\mathcal X}}
\newcommand{\cT}{\ensuremath{\mathcal T}}

\newcommand{\R}{\ensuremath{\mathbb R}}
\newcommand{\C}{\ensuremath{\mathbb C}}

\newcommand{\PP}{\ensuremath{\mathbb P}}
\newcommand{\Z}{\ensuremath{\mathbb Z}}
\newcommand{\Q}{\ensuremath{\mathbb Q}}

\newcommand{\kahler}{K\"ahler\xspace}
\newcommand{\hk}{hyperk\"ahler\xspace}

\newcommand{\qk}{quaternionic-K\"ahler\xspace}

\newcommand{\I}{{\mathrm i}}

\newcommand{\de}{\mathrm{d}}

\newcommand{\gauge}{{\mathrm{g}}}
\newcommand{\flav}{{\mathrm{f}}}

\newcommand{\abs}[1]{\lvert#1\rvert}
\newcommand{\norm}[1]{\lVert#1\rVert}
\newcommand{\IP}[1]{\langle#1\rangle}

\newcommand{\ti}[1]{\textit{#1}}

\newcommand{\In}{{\mathrm{in}}}
\newcommand{\Out}{{\mathrm{out}}}

\DeclareMathOperator{\re}{Re}
\DeclareMathOperator{\Sp}{Sp}

\DeclareMathOperator{\Hom}{Hom}
\DeclareMathOperator{\TChar}{TChar}

\DeclareMathOperator{\Stab}{Stab}

\renewcommand{\sf}{\mathrm{sf}}

\newcommand{\insfig}[3]{\begin{figure}[htbp] \centering \includegraphics[width=#1]{figures/#2-crop.pdf} \caption{#3} \label{fig:#2} \end{figure}}

\newenvironment{vardescription}[1]%
  {\begin{list}{}{%
    \settowidth{\labelwidth}{\textbf{#1:}}%
    \setlength{\leftmargin}{\labelwidth}\addtolength{\leftmargin}{\labelsep}}}%
  {\end{list}}

\begin{document}

\bibliographystyle{utphys}
\setcounter{page}{1}

\title{Notes on a new construction of \hk metrics}
\date{}
\author{Andrew Neitzke \\ \\ \small Department of Mathematics \\ \small University of Texas at Austin}

\maketitle

\section{Overview}

In joint work with Davide Gaiotto and Greg Moore \cite{Gaiotto:2008cd} we recently proposed
a new connection between \hk geometry and the counting
of BPS states in supersymmetric field theory.  While the story is
motivated by physics, it leads to a concrete new recipe for constructing complete
\hk metrics on the total spaces of certain complex integrable systems.

The aim of this note is briefly to describe what this recipe is, and to
comment on some of the issues involved in converting it
into an actual theorem.

\medskip

Let us briefly describe some of the highlights.

\begin{itemize}
 \item We begin with a collection of ``integrable system data'' described in Section \ref{sec:int-data} below.
These data include a complex manifold
$\cB$ containing a divisor $D$.  For example, $\cB$ could be the complex plane, and $D$ some collection
of points.  The data also include a local system of lattices $\Gamma$ over $\cB' = \cB \setminus D$, from which we build
a $2r$-torus bundle $\cM'$ over $\cB'$, with nontrivial monodromy
around $D$.  Finally, we have a ``central charge'' homomorphism $Z: \Gamma \to \C$, varying holomorphically over $\cB'$.
From these data we build a simple explicit \hk metric $g^\sf$ on $\cM'$.
However, the metric $g^\sf$ is incomplete, and our main interest is in complete metrics.

\item Naively we might hope to complete $g^\sf$ by adding some degenerate torus fibers
over $D$, thus extending $\cM'$ to $\cM \supset \cM'$,
in such a way that $g^\sf$ will extend to $\cM$.  However, it seems that this is impossible:  roughly
speaking, $g^\sf$ is too homogeneous to have such an extension.  Instead, we construct a new metric $g$ on $\cM'$, which differs from $g^\sf$ by certain ``quantum corrections.''

\item The quantum corrections are obtained
by solving a certain explicit integral equation, \eqref{inteq} below.
The main new ingredient in this equation is
a set of integer ``invariants'' $\Omega(\gamma)$, which should be
examples of generalized Donaldson-Thomas invariants in the sense of \cite{ks1,Joyce:2008pc}.
In particular, the Kontsevich-Soibelman wall-crossing formula for generalized Donaldson-Thomas invariants,
as written in \cite{ks1}, plays an important role in the construction.
Indeed the original motivation for this construction was an attempt to understand
the physical meaning of the formula of \cite{ks1}.

\item  Both the metrics $g^\sf$ and $g$ depend on a real parameter $R > 0$;
in the limit as $R \to \infty$, the torus fibers of $\cM'$ collapse, in either metric.
The corrections $g - g^\sf$ are
exponentially suppressed in $R$ when we are away from $D$:  so as $R \to \infty$,
$g$ looks very close to $g^\sf$ except in a small neighborhood of the singular fibers.
Near the singular fibers the quantum corrections become large, and in particular we
expect that the corrected $g$ can be extended over the singular fibers.

This description of $g$ near the $R \to \infty$ limit
should be thought of as an example of a more general picture of the
geometry of Calabi-Yau manifolds near their large complex structure limit,
proposed by Gross-Wilson \cite{MR1863732},
Kontsevich-Soibelman \cite{MR1882331} and Todorov,
motivated by the Strominger-Yau-Zaslow picture of mirror symmetry \cite{Strominger:1996it}.

\item In many examples where our recipe
can be applied, it turns out that the \hk metrics in question were already known to exist.
The example we have studied in most detail is
that of rank-2 Hitchin systems with semisimple ramification \cite{Gaiotto:2009hg}.
We briefly describe that example in Section \ref{sec:hitchin} below.

\item Our recipe has not really been tested so far, in the sense that nobody has tried hard to
use it to get new explicit information about interesting \hk metrics.
We believe that this should be possible:  at the very least it should be possible to get
a precise asymptotic series for $g$ as $R \to \infty$.

\end{itemize}

\medskip

I thank Cesar Garza, Tom Sutherland, and the anonymous referee for very useful comments and for
correcting several errors in an earlier draft of this note.  I would also like to thank
Davide Gaiotto and Greg Moore for a very enjoyable collaboration.

This work is supported by National Science Foundation grants DMS-1160461 and DMS-1151693.

\section{Integrable system data}

\subsection{Data} \label{sec:int-data}

Our construction begins with the following data:

\begin{vardescription}{Data 3c}
\begin{shaded}
\item[Data 1] A complex manifold $\cB$, of dimension $r$ (``Coulomb branch'').
\end{shaded}
\begin{shaded}
\item[Data 2] A divisor $D \subset \cB$ (``discriminant locus'').
Let $\cB' = \cB \setminus D$ (``smooth locus'').  We use $u$ to denote a general point of $\cB'$.
\end{shaded}
\begin{shaded}
\item[Data 3a] A local system $\Gamma_\gauge$ over $\cB'$, with fiber a rank-$2r$ lattice, equipped with
a nondegenerate antisymmetric integer-valued pairing $\IP{,}$.
Abusing notation we will also
use $\IP{,}$ to denote the inverse pairing on $\Gamma_\gauge^*$ (not necessarily integer-valued.)
\item[Data 3b] A fixed lattice $\Gamma_\flav$ (possibly trivial).  We sometimes think of $\Gamma_\flav$
as the fiber of a trivial local system of lattices over $\cB'$.
\item[Data 3c] A local system $\Gamma$ of lattices over $\cB'$, given as an extension
\begin{equation} \label{latticeext}
 0 \to \Gamma_\flav \to \Gamma \to \Gamma_\gauge \to 0.
\end{equation}
The pairing $\IP{,}$ on $\Gamma_\gauge$ induces one on $\Gamma$ which we also denote $\IP{,}$.
The radical of this pairing is $\Gamma_\flav$.
\end{shaded}
\begin{shaded}
\item[Data 4] A homomorphism $Z: \Gamma \to \C$, varying holomorphically over $\cB'$.  For any
local section $\gamma$ of $\Gamma$ we thus get a local holomorphic function $Z_{\gamma}$ on $\cB'$.
\end{shaded}

\begin{shaded}
\item[Data 5] A homomorphism $\theta_\flav: \Gamma_\flav \to \R / 2 \pi \Z$.
\end{shaded}
\end{vardescription}

These data are subject to several conditions:
\begin{vardescription}{Condition 1}
\begin{shaded}
\item[Condition 1] $Z_{\gamma_f}$ is a constant function on $\cB'$ for any $\gamma_f \in \Gamma_\flav$.  (As a consequence, the $\Gamma^*$-valued 1-form $\de Z$ actually descends to $\Gamma_\gauge^*$; we use this
in formulating Condition 2.)
\end{shaded}
\begin{shaded}
\item[Condition 2] $\IP{\de Z \wedge \de Z} = 0$.
\end{shaded}
\begin{shaded}
\item[Condition 3] For any $u \in \cB'$, the $\de Z_\gamma(u)$ span $\cT^*_u \cB'$.
\end{shaded}
\end{vardescription}

\subsection{Integrable system}

The above data are enough to determine an incomplete complex integrable system, i.e. a holomorphic
symplectic manifold $\cM'$ which is a fibration over a complex base manifold $\cB'$, with fibers
complex Lagrangian tori.  We now describe $\cM'$.

For any fiber $\Gamma_u$ of $\Gamma$,
let $\TChar_u(\Gamma, \theta_\flav)$ be the set of twisted unitary characters of $\Gamma_u$, i.e. maps
$\theta: \Gamma_u \to \R / 2 \pi \Z$ obeying
\begin{equation} \label{twisted-character}
 \theta_\gamma + \theta_{\gamma'} = \theta_{\gamma + \gamma'} + \pi \IP{\gamma, \gamma'},
\end{equation}
agreeing with $\theta_\flav$ when restricted to $\Gamma_\flav \subset \Gamma_u$.
$\TChar_u(\Gamma, \theta_\flav)$ is topologically a torus $(S^1)^{2r}$.  Letting $u$ vary, the
$\TChar_u(\Gamma, \theta_\flav)$ are the fibers of a torus bundle $\cM'$ over $\cB'$.
Any local section $\gamma$ of $\Gamma$ then gives a function ``evaluation on $\gamma$,''
\begin{equation}
 \theta_\gamma: \cM' \to \R / 2 \pi \Z.
\end{equation}
These are the angular coordinates on the torus fibers of $\cM'$.

Now we want to construct the complex structure and
holomorphic symplectic form on $\cM'$.  For this purpose note that
we have canonical functions
\begin{equation}
 Z_\gamma: \cM' \to \C,
\end{equation}
pulled back from the base $\cB'$.
Differentiating gives a collection of 1-forms $\de \theta_\gamma$ and $\de Z_\gamma$ on $\cM'$, which are linear in $\gamma$
and vanish for $\gamma \in \Gamma_\flav$, and hence can be organized into
$\Gamma_\gauge^*$-valued 1-forms $\de\theta$ and $\de Z$.
Then define the complex 2-form
\begin{equation} \label{omegaplus}
 \omega_+ = - \frac{1}{2\pi} \IP{\de Z \wedge \de\theta}.
\end{equation}
There is a unique complex structure on $\cM'$ for which $\omega_+$ is of type $(2,0)$.
We call this complex structure $J(\zeta = 0)$, for a reason which will emerge momentarily.

The two-form $\omega_+$ gives a holomorphic symplectic structure on $(\cM', J(\zeta = 0))$.  With respect to this structure,
the projection $\pi: \cM' \to \cB'$ is holomorphic, and
the torus fibers $\cM'_u = \pi^{-1}(u)$ are compact complex Lagrangian submanifolds.

\subsection{Affine structures}

Although we do not use it explicitly in the rest of this note,
it may be useful to mention that our data
determine an $S^1$ worth of (symplectic) affine structures on $\cB'$.
Fix some $\vartheta \in \R / 2 \pi \Z$.  Then pick a patch $U \subset \cB'$
on which $\Gamma_\gauge$ admits a basis of local sections,
$\gamma_1, \dots, \gamma_{2r}$, in which $\IP{,}$ is the standard symplectic pairing.
Also choose a local splitting $\rho: \Gamma_\gauge \to \Gamma$ of \eqref{latticeext}.
Then the functions
\begin{equation}
 f_i = \re(e^{\I \vartheta} Z_{\rho(\gamma_i)})
\end{equation}
are local coordinates on $U$ (possibly after shrinking $U$).
The transition functions on overlaps $U \cap U'$
are valued in $\Sp(2r, \Z) \ltimes \R^{2r}$
(the $\Sp(2r, \Z)$ part comes from the choice of basis of $\Gamma_\gauge$,
the $\R^{2r}$ from the choice of splitting $\rho$.)

\section{Semiflat hyperk\"ahler metric}

We now impose one more condition.  Recall that a positive 2-form $\omega$ on a complex manifold
is a real 2-form for which $\omega(v, Jv) > 0$ for all real tangent vectors $v$.

\begin{vardescription}{Condition 4}
\begin{shaded}
\item[Condition 4] $\IP{\de Z \wedge \de\bar{Z}}$ is a positive 2-form on $\cB'$.
\end{shaded}
\end{vardescription}

\subsection{Semiflat metric}

Fix $R \in \R_+$.
$\cM'$ carries a canonical 2-form,
\begin{equation} \label{omega3sf}
\omega_3^\sf = \frac{R}{4} \langle \de Z \wedge \de \bar{Z} \rangle - \frac{1}{8 \pi^2 R} \langle \de\theta \wedge \de\theta \rangle.
\end{equation}
This form is of type $(1,1)$ in complex structure $J(\zeta = 0)$ and positive.  So the triple
$(\cM', J(\zeta = 0), \omega_3^\sf)$ determine a \kahler metric $g^\sf$ on $\cM'$.
In fact this metric is \hk.  As far as I know, the first place where this was
shown is in \cite{Cecotti:1989qn} (albeit in somewhat different notation); see also \cite{Freed:1997dp}
for a more modern account.  Alternatively, though, the \hk property is a consequence
of the twistorial construction of the metric which we will give below.

The superscript $^\sf$ stands for ``semi-flat'':  this terminology first appeared in \cite{MR1863732},
where it was used to refer to an important special case introduced in \cite{Greene:1990ya}.
The reason for the name is that
$g^\sf$ is flat when restricted to any torus fiber $\cM'_u$, and $\cM'_u$ has half the dimension of $\cM'$.

\subsection{Twistorial description of the semiflat metric}

Let us now describe a different, ``twistorial'' way of constructing the metric $g^\sf$;
this alternative description is what we will generalize in
our construction of the quantum-corrected metric $g$ below.

Any \hk metric on a manifold $\cM'$ determines --- and is determined by ---
a collection of holomorphic symplectic structures $(\cM', J(\zeta), \varpi(\zeta))$
labeled by $\zeta \in \C\PP^1$.  In the general theory of \hk manifolds all $\zeta$ are
on the same footing.  However, for the \hk manifolds we are describing in this note,
the points $\zeta = 0$ and $\zeta = \infty$ will play
a distinguished role.  It is then convenient to expand $\varpi(\zeta)$ as
\begin{equation} \label{eq:Omegazeta}
\varpi(\zeta) = - \frac{i}{2\zeta} \omega_+ + \omega_3 - \frac{i}{2} \zeta \overline{\omega_+}
\end{equation}
where $\omega_+$, $\omega_3$ are respectively the holomorphic symplectic form
and \kahler form, both relative to the complex structure $J(\zeta = 0)$.

In the particular case of the \hk metric $g^\sf$, we have written these 2-forms
explicitly above in \eqref{omegaplus}, \eqref{omega3sf}.
We will now give an alternative description of the holomorphic symplectic forms
$\varpi(\zeta)$ corresponding to $g^\sf$, roughly by exhibiting explicit
holomorphic Darboux coordinates.

Let $T_u$ denote the complex torus of twisted complex
characters of $\Gamma_u$.  $T_u$ has canonical $\C^\times$-valued
functions $X_\gamma$ ($\gamma \in \Gamma_u$)
obeying
\begin{equation}
 X_\gamma X_{\gamma'} = (-1)^{\IP{\gamma,\gamma'}} X_{\gamma + \gamma'},
\end{equation}
and a Poisson structure
\begin{equation}
 \{X_\gamma, X_{\gamma'}\} = \IP{\gamma, \gamma'} X_{\gamma+\gamma'}.
\end{equation}
The $T_u$ glue together into a local system over $\cB'$ with
fiber a complex Poisson torus.  Let $T$ denote the pullback of this local system to $\cM'$.

Now we consider a section $\cX^\sf$ of $T$, depending on
an auxiliary parameter $\zeta \in \C^\times$.
Locally this just means a collection of ``coordinate'' functions
\begin{equation}
 \cX^\sf_\gamma: \cM' \times \C^\times \to \C^\times
\end{equation}
(defined by $\cX^\sf_\gamma = (\cX^\sf)^* X_\gamma$, with $\gamma$ a local section of $\Gamma$).
We often write these functions
as $\cX^\sf_\gamma(\zeta)$, leaving the $\cM'$ dependence implicit.
$\cX^\sf(\zeta)$ is given by a simple closed formula:
\begin{equation} \label{eq:xsf}
 \cX^\sf_\gamma(\zeta) = \exp \left[ \pi R \frac{Z_\gamma}{\zeta} + \I\theta_\gamma + \pi R \zeta \bar{Z}_\gamma \right].
\end{equation}
Now a direct computation shows\footnote{Note that this computation uses Condition 2,
the fact $\IP{\de Z \wedge \de Z} = 0$ --- if we did not impose this condition, then
computing the right side of \eqref{eq:hsympsf} would produce a term
$\IP{\de Z \wedge \de Z} / \zeta^2$, which would not match the form
of $\varpi^\sf(\zeta)$.}
\begin{equation} \label{eq:hsympsf}
 \varpi^\sf(\zeta) = \frac{1}{8 \pi^2 R} \IP{\de \log \cX^\sf(\zeta) \wedge \de \log \cX^\sf(\zeta)}.
\end{equation}
So the $\cX^\sf_\gamma(\zeta)$ are ``holomorphic Darboux coordinates'' on $\cM'$,
determining the holomorphic symplectic structures for all $\zeta \in \C^\times$,
and hence the \hk metric $g^\sf$.  In short:
\ti{knowing the functions $\cX^\sf_\gamma(\zeta)$ is equivalent to knowing the \hk metric $g^\sf$.}

A global way of thinking about this construction is to say that for each $\zeta \in \C^\times$
we \ti{pull back} the structure
of holomorphic Poisson manifold from $T$ to $\cM'$, using the section $\cX^\sf(\zeta)$ of $T$.\footnote{Of course $T$ is a local system of tori, not a single torus, so the last sentence does not strictly make
sense; but locally we can view $\cX^\sf(\zeta)$ as a map into the space of local flat sections of $T$,
which \ti{is} a single holomorphic Poisson torus.}
After pullback the Poisson structure is actually nondegenerate, i.e. it arises from a holomorphic
symplectic structure.

\section{Instanton corrections to $\cX$}

We explained above how the semiflat section $\cX^\sf$ can be used to construct the holomorphic-symplectic
forms $\varpi(\zeta)$ corresponding to the semiflat metric $g^\sf$.
We now want to construct a new, ``quantum-corrected'' section $\cX$.  In the next section we will
use $\cX$ to build a quantum-corrected metric $g$.

\subsection{BPS degeneracies and Riemann-Hilbert problem}

The key new ingredient determining the quantum corrections is:

\begin{vardescription}{Data 6}
\begin{shaded}
 \item[Data 6] A function $\Omega: \Gamma \to \Z$.
\end{shaded}
\end{vardescription}

For each local section $\gamma$ of $\Gamma$ this gives a locally defined integer-valued function
$\Omega(\gamma)$ on $\cB'$.  I emphasize that $\Omega$ is \ti{not} required to be continuous:
indeed Condition 7 below will imply that it is generally not continuous, but jumps in a specific
way (governed by the Kontsevich-Soibelman
wall-crossing formula) at real-codimension-1 loci in $\cB'$.

$\Omega$ should obey a simple parity-invariance condition:
\begin{vardescription}{Condition 5}
\begin{shaded}
 \item[Condition 5] $\Omega(\gamma;u) = \Omega(-\gamma;u)$.
\end{shaded}
\end{vardescription}

We can now formulate the key ingredient in our construction, a certain Riemann-Hilbert problem.
We need a little notation.
Any $\gamma \in \Gamma_u$ gives a birational Poisson automorphism $\cK_\gamma$ of
$T_u$, defined by
\begin{equation}
 \cK_\gamma^* X_{\gamma'} = X_{\gamma'} (1 - X_\gamma)^{\IP{\gamma,\gamma'}}.
\end{equation}
$\cK_\gamma$ and $\cK_{\gamma'}$ commute if and only if $\IP{\gamma, \gamma'} = 0$.
Define a ray associated to each $\gamma \in \Gamma_u$,
\begin{equation}
 \ell_\gamma(u) := Z_\gamma(u) \R_-.
\end{equation}
Then to each ray $\ell$ running from the origin to infinity in the $\zeta$-plane,
associate a certain birational Poisson automorphism of $T_u$ (first written down in \cite{ks1}),
\begin{equation}
S_\ell(u) := \prod_{\gamma: \ell_\gamma(u) = \ell} \cK_{\gamma}^{\Omega(\gamma;u)}.
\end{equation}
We call the $\ell$ for which $S_\ell(u) \neq 1$ ``BPS rays.''
Finally, we define an antiholomorphic involution $\rho$ of $T_u$ by
\begin{equation}
 \rho^* X_\gamma = \overline X_{-\gamma}.
\end{equation}

Now we can formulate the Riemann-Hilbert problem.  Fix $u \in \cB'$.
We seek a map
\begin{equation}
 \cX: \cM_u \times \C^\times \to T_u
\end{equation}
with the following properties:

\begin{enumerate}

\item $\cX$ depends piecewise-holomorphically on $\zeta \in \C^\times$, with discontinuities only at the
rays $\ell_\gamma(u)$ for $\gamma \in \Gamma_u$ with $\Omega(\gamma;u) \neq 0$.

\item
The limits $\cX^\pm$ of $\cX$ as $\zeta$ approaches any ray $\ell$ from both sides exist and are related by
\begin{equation} \label{x-jumps}
 \cX^+ = S_\ell^{-1} \circ \cX^-.
\end{equation}

\item
$\cX$ obeys the reality condition
\begin{equation} \label{reality-condition}
 \cX(-1/\bar\zeta) = \rho^* \cX(\zeta).
\end{equation}

\item
For any $\gamma$, $\lim_{\zeta \to 0} \cX_\gamma(\zeta) / \cX^\sf_\gamma(\zeta)$ exists and is real.

\suspend{enumerate}

We expect that the $\cX$ with these properties should be unique if it exists, by analogy with what
is known for similar Riemann-Hilbert problems appearing in \cite{MR1213301,Cecotti:1993rm}.

\subsection{Solving the Riemann-Hilbert problem}

To find a solution of these conditions we contemplate the integral equation
\begin{equation} \label{inteq}
\cX_\gamma(x, \zeta) = \cX^\sf_\gamma(x, \zeta) \exp \left[ -\frac{1}{4
\pi \I} \sum_{\gamma'} \Omega(\gamma';u) \langle \gamma,\gamma'
\rangle \int_{\ell_{\gamma'}(u)} \frac{d\zeta'}{\zeta'} \frac{\zeta' +
\zeta}{\zeta' - \zeta} \log (1-\cX_{\gamma'}(x, \zeta'))\right].
\end{equation}
For any fixed $x \in \cM'$, \eqref{inteq} is a functional equation for the functions
$\cX_\gamma(x, \cdot): \C^\times \to \C^\times$.
We claim that if we find a collection of functions $\cX_\gamma$ obeying this equation, they are a solution
of our Riemann-Hilbert problem (in other words they obey the 4 conditions set out in the
last section).

A natural way to try to
produce a solution of \eqref{inteq} is by iteration, beginning with $\cX = \cX^\sf$.
In \cite{Gaiotto:2008cd} we sketch a proof that this iteration indeed converges for large enough $R$, to
the \ti{unique} solution of \eqref{inteq},
under ``reasonable'' growth conditions on the $\Omega(\gamma;u)$ (stated more precisely in \cite{Gaiotto:2008cd}):

\begin{vardescription}{Condition 6}
\begin{shaded}
 \item[Condition 6] $\Omega(\gamma;u)$ does not grow too quickly as a function of $\gamma$ for
fixed $u$.
\end{shaded}
\end{vardescription}

Let us make a few remarks about this:
\begin{itemize}
\item
This approach to the Riemann-Hilbert problem was inspired by the treatment of a similar problem in
\cite{MR1213301,Cecotti:1993rm}.
At least morally speaking, ours is an infinite-dimensional version of the one discussed there, with the group
$GL(K,\R)$ replaced by the group of symplectomorphisms of the torus $T$.

\item Our arguments are not strong enough to give \ti{uniform} convergence of the iteration as we
vary $u$, since $\Omega(\gamma;u)$ and $Z_\gamma(u)$ depend on $u$; in particular, the correct notion of ``large enough $R$'' may depend on $u$.  Roughly speaking, the speed of the convergence is set by the largest $e^{- 2 \pi R \abs{Z_\gamma(u)}}$ for which $\Omega(\gamma;u) \neq 0$.

\item We did not give a complete proof that the $\cX_\gamma$ obey our asymptotic Condition 4; we expect though
that it should be possible to prove it directly, at least for large enough values of the parameter $R$,
along similar lines to what was discussed in \cite{MR1213301,Cecotti:1993rm}.
Essentially the idea is that for large $R$ the integrals in \eqref{inteq}
have a finite limit as $\zeta \to 0$:  this is easy to check directly if we replace $\cX$ by $\cX^\sf$, and we expect
that this property should be preserved by the iteration.

\item The $\cX_\gamma$ are ``quantum-corrected'' versions of the
original functions $\cX^\sf_\gamma$.
As with the $\cX^\sf_\gamma$, the $\cX_\gamma$ can be thought of
as $\cX_\gamma = \cX^* X_\gamma$ for some section $\cX$ of the complex torus bundle $T$.

\end{itemize}

\subsubsection{Sums over trees} \label{sec:trees}

We also give a formula for a solution $\cX_\gamma$ of \eqref{inteq} as a sum over certain iterated integrals,
as follows.  (It is not clear at the moment whether this sum actually
\ti{converges} or gives instead only an asymptotic series.)

We first introduce $\Q$-valued invariants related to the $\Omega(\gamma)$ by a ``multi-cover formula'' \cite{ks1},
\begin{equation} \label{eq:dt-inv}
c(\gamma) = \sum_{n = 1}^\infty \frac{\Omega(\gamma/n)}{n^2}.
\end{equation}
(Here we take $\Omega(\gamma/n) = 0$ by definition whenever $n$ does not divide $\gamma$.)
We consider rooted trees, with
edges labeled by pairs $(i,j)$ (where $i$ is the node closer to the root),
and each node decorated by some $\gamma_i \in \Gamma$.
Let $\cT$ denote such a tree.
Define a weight attached to $\cT$ by
\begin{equation}
c(\cT) = \frac{1}{\abs{\mathrm{Aut}(\cT)}} \prod_{i \in {\mathrm{Nodes}}(\cT)} c(\gamma_i) \prod_{(i,j) \in {\mathrm {Edges}}(\cT)} \IP{\gamma_i, \gamma_j}.
\end{equation}
Let $\gamma_\cT$ denote the decoration at the root node of $\cT$.
We define a function $\cG_\cT(x,\zeta)$ on (a patch of) $\cM$ inductively as follows:
deleting the root node from $\cT$ leaves behind a set of trees $\cT_a$, and
\begin{equation}
\cG_\cT(x,\zeta) = \frac{1}{4 \pi \I} \int_{\ell_{\gamma_\cT}} \frac{\de \zeta'}{\zeta'} \frac{\zeta' + \zeta}{\zeta' - \zeta} \cX^\sf_{\gamma_\cT}(x,\zeta') \prod_a \cG_{\cT_a}(x,\zeta').
\end{equation}
Then a formal solution of \eqref{inteq} can be given as
\begin{equation} \label{eq:iterative-X}
\cX_\gamma(x,\zeta) = \cX^\sf_\gamma (x,\zeta)\exp \left[ \sum_\cT
\langle \gamma, \gamma_\cT \rangle c(\cT) \cG_\cT(x,\zeta) \right].
\end{equation}

\subsection{Wall-crossing formula}

Define the ``locus of marginal stability'' by
\begin{equation}
 W = \{ u: \exists \gamma_1, \gamma_2 \text{ with } \Omega(\gamma_1;u) \neq 0, \Omega(\gamma_2;u) \neq 0, Z_{\gamma_1}(u) / Z_{\gamma_2}(u) \in \R_+ \} \subset \cB'.
\end{equation}
This $W$ is a union of countably many components (``walls'') each of which has real codimension $1$ in $\cB'$.  For our construction to work,
the integers $\Omega(\gamma;u)$ must \ti{jump} as $u$ crosses any of these walls.
More precisely, they must jump in accordance with the celebrated
wall-crossing formula of Kontsevich and Soibelman
\cite{ks1}.  We now describe this formula, essentially following \cite{ks1}, with a few slight adaptations
to our context.

Let $V$ be a strictly
convex cone in $\C$ with apex at the origin.  Then for any $u \notin W$, define
\begin{equation} \label{av}
 A_V(u) = \prod_{\gamma: Z_\gamma(u) \in V} \cK_\gamma^{\Omega(\gamma; u)} = \prod_{\ell \subset V} S_\ell(u),
\end{equation}
where the product is taken in order of increasing $\arg Z_\gamma(u)$.
$A_V(u)$ is a birational
Poisson automorphism of $T_u$.\footnote{This statement needs a little amplification
since the product in \eqref{av} may be infinite.  One should more precisely
think of $A_V(u)$ as living in a
certain prounipotent completion of the group generated by $\{ \cK_\gamma \}_{\gamma: Z_\gamma(u) \in V}$
as explained in \cite{ks1}.}
Knowing $A_V(u)$ is sufficient to determine the
$\Omega(\gamma; u)$ for $\gamma$ with $Z_\gamma(u) \in V$; thus
we can think of $A_V(u)$ as a sophisticated kind of generating function.

Define a \ti{$V$--good path} to be a path $p \subset \cB'$
along which there is no point $u$ with $Z_\gamma(u) \in \partial V$ and $\Omega(\gamma;u) \neq 0$.
(So as we travel along a $V$--good path, no BPS rays enter or exit $V$.)

\begin{vardescription}{Condition 7}
\begin{shaded}
 \item[Condition 7] If $u$ and $u'$ are the endpoints of a $V$--good path $p$,
then $A_V(u)$ and $A_V(u')$ are related by parallel transport in $T$ along $p$.
\end{shaded}
\end{vardescription}

Condition 7 is essentially the wall-crossing formula of Kontsevich and Soibelman \cite{ks1}.  It
is strong enough to determine all $\Omega(\gamma; u)$, if we have Data 1--4 and also know the
$\Omega(\gamma; u_0)$ for some fixed $u_0$.
In fact, at first sight it might seem to imply simply that $\Omega(\gamma; u)$ are locally constant functions
of $u$ on $\cB'$.  This is almost right:  what it actually implies is that $\Omega(\gamma; u)$ are locally constant functions of $u$
on $\cB' \setminus W$.
The point is that when $u$ hits $W$ the order of the factors in the product \eqref{av} is changed; as a result, for $A_V$ to remain
constant, the individual factors must in general also change.  In other words, the $\Omega(\gamma;u)$ must jump.

Condition 7 determines precisely how the $\Omega(\gamma; u)$ jump
when $u$ crosses some component of $W$.  It is in this sense
that it is a wall-crossing formula.

\subsection{Absence of unwanted jumps in $\cX$} \label{sec:no-jumps}

Under this condition, let us revisit the solution $\cX$ of the
Riemann-Hilbert problem, and now vary the point $u \in \cB$
as well as $\zeta \in \C^\times$.  We have already noted that for any fixed $u$, $\cX$ is discontinuous along
the BPS rays.  Letting $u$ vary this becomes the statement that $\cX$ is discontinuous along the locus
\begin{equation}
 L = \left\{ (u, \zeta): \exists \gamma \in \Gamma_u \textrm{ with } Z_{\gamma}(u) / \zeta \in \R_- \textrm{ and } \Omega(\gamma;u) \neq 0 \right\} \subset \cB' \times \C^\times.
\end{equation}
If Condition 7 is not obeyed, it is straightforward to show that these cannot be the only discontinuities of $\cX$:
there must be additional jumps when $u$ meets the walls of marginal stability $W \subset \cB'$.
Such additional jumps would be a problem for our construction of the corrected \hk metric below.

On the other hand, if Condition 7 \ti{is} obeyed, then we claimed in \cite{Gaiotto:2008cd} that $\cX$ is actually continuous.
This statement would follow directly from uniqueness of the solution of our Riemann-Hilbert problem,
since Condition 7 says that
the two Riemann-Hilbert problems we obtain by approaching the wall $W$ from two sides are actually the same.

\section{Corrected metric}

\subsection{Construction}

Having defined the section $\cX(\zeta)$ of $T$, we are ready to describe the corrected \hk metric $g$.
The idea is similar to one we used above in our description of $g^\sf$.  Namely, for each $\zeta \in \C^\times$,
we use $\cX(\zeta)$ to pull back a holomorphic symplectic structure $\varpi^{(\zeta)}$ from $T$ to $\cM'$.
As we have noted, $\cX(\zeta)$ is not continuous; it has jumps along the locus
$\pi^{-1}(L) \subset \cM' \times \C^\times$, given by \eqref{x-jumps}.  Fortunately this jump is by composition
with a Poisson morphism of $T$, and thus does not affect $\varpi^{(\zeta)}$.  So $\varpi^{(\zeta)}$
is continuous, and depends holomorphically on $\zeta \in \C^\times$.
In order to define an honest holomorphic symplectic structure, $\varpi^{(\zeta)}$ should also be
nondegenerate.  One expects this to be true at least
for large enough $R$, since it is true for $\cX^\sf$ and
$\cX$ differs from $\cX^\sf$ only by corrections that are exponentially suppressed
at large $R$.

Now our key claim is that

\medskip

\ti{$\varpi^{(\zeta)}$ is of the form \eqref{eq:Omegazeta},
where $(\omega_\pm, \omega_3)$ are symplectic forms defining a \hk structure on
$\cM'$.}

\medskip

This is our construction of the new \hk metric $g$ on $\cM'$.

\subsection{Twistor space}

Let us make a few comments about how the claim above is motivated.
One obvious necessary condition is $\varpi^{(-1/\bar\zeta)} = \overline{\varpi^{(\zeta)}}$.
This follows from the reality condition \eqref{reality-condition} (property 3 of the Riemann-Hilbert problem).
We also need to see that $\varpi^{(\zeta)}$ has only
a simple pole at $\zeta = 0$ (hence also at $\zeta = \infty$.)  This follows
from our asymptotic condition on $\cX$ (property 4 of the Riemann-Hilbert problem).
So $\varpi^{(\zeta)}$ indeed determines a complex 2-form
$\omega_+$ and a real 2-form $\omega_3$.
Of course this is still not enough to guarantee that these 2-forms fit together into an \hk structure on $\cM'$.
The most delicate point is to show that indeed they do.

For this we use the ``twistor space'' construction
\cite{Hitchin:1986ea,MR1206066}.  We consider the space $\cZ = \cM \times \C\PP^1$.  The 2-form
$\varpi$ equips $\cZ$ with a complex structure for which the projection to $\C\PP^1$ is holomorphic,
and a fiberwise holomorphic symplectic form (globally twisted by $\cO(2)$),
obeying an appropriate reality condition.  Moreover $\cZ$ has a family of distinguished holomorphic sections
labeled by points $x \in \cM'$, given by the tautological-looking formula $s_x(\zeta) = (x, \zeta)$.
In this situation, the twistor space construction promises us a \hk metric on $\cM'$,
provided that the normal bundle $N(s_x)$ to each such
section is a direct sum of copies of $\cO(1)$.  This condition on the normal bundle is the most delicate
part of the story; we argue in \cite{Gaiotto:2008cd} that it is a consequence of the asymptotic conditions obeyed by the
section $\cX$ as $\zeta \to 0$.

\subsection{Improvement of singularities}

So far we have described how to construct a ``quantum-corrected \hk metric'' $g$ on $\cM'$.
The reader may be wondering why we have bothered to do so much work.  After all, we already had a
perfectly good \hk metric $g^\sf$ on $\cM'$.

However, $g^\sf$ has one important deficiency (in all but
the most trivial examples):  it is incomplete.  The reason for this incompleteness is the fact that $g^\sf$
is defined only on $\cM'$, which has smooth torus fibers over all points of $\cB' \subset \cB$, but does not
include fibers over points of the ``singular locus'' $D \subset \cB$.
Typically one can complete $\cM'$ topologically to a natural $\cM$, with a projection $\pi: \cM \to \cB$,
such that the fiber over a point of $D$ is some kind of degenerate torus.
One might then try to extend $g^\sf$ to a metric on the whole $\cM$.
This however appears to be impossible.

One answer to the question ``why is $g$ better than $g^\sf$?''
is that, if $\Omega$ is chosen appropriately, we expect that
$g$ does admit an extension to a metric on $\cM$, which in many cases will be complete.
So morally the statement is that the quantum corrections ``improve'' the behavior of the metric
near the singular locus $D$.  We will discuss an example in the next section.

\section{Ooguri-Vafa metric}

In \cite{Gaiotto:2008cd} we discussed a model example of this phenomenon of improvement of singularities.
Fix some constant $\Lambda \in \C$ (which enters the story in a trivial way:  it is safe to fix $\Lambda = 1$
if you prefer.)
We choose our data as follows:
\begin{vardescription}{Data 6}
\item[Data 1] $\cB$ is the disc $\{ \abs{u} < \abs{\Lambda} \}$.
\item[Data 2] The discriminant locus is $D = \{u = 0\} \subset \cB$.  So $\cB'$ is the punctured disc.
\item[Data 3] $\Gamma = \Gamma_\gauge$ is a local system of rank-2 lattices over $\cB'$.  With respect to a local
basis of sections $(\gamma_m, \gamma_e)$, with $\IP{\gamma_m, \gamma_e} = 1$, the monodromy around the puncture $u = 0$
is $\gamma_e \to \gamma_e, \gamma_m \to \gamma_m + \gamma_e$.  $\Gamma_\flav$ is trivial.
\item[Data 4] With respect to the same local basis of sections,
$Z_{\gamma_e}(u) = u$, $Z_{\gamma_m}(u) = \frac{1}{2 \pi i} (u \log \frac{u}{\Lambda} - u)$.  Note that analytically
continuing around $u = 0$ we get $Z_{\gamma_m} \to Z_{\gamma_m} + Z_{\gamma_e}$, consistent with the monodromy of $\Gamma$;
in other words $Z$ is really globally defined.
\item[Data 5] Since $\Gamma_\flav$ is trivial, $\theta_\flav$ is trivial.
\item[Data 6] For all $u$, we have $\Omega(\gamma;u) = \begin{cases} 1 \text{ for } \gamma \in \{ \gamma_e, -\gamma_e \}, \\ 0 \text{ otherwise.} \end{cases}$
\end{vardescription}

In this case our construction can be carried out very explicitly (for any value of the parameter $R$):
the integral equation \eqref{inteq} becomes simply an integral \ti{formula}, or said otherwise, the iterative procedure of
finding a solution actually terminates after a single step.  So in this case we know the functions $\cX_\gamma$ exactly.
Applying our construction then yields an \hk metric $g$ on a torus fibration $\cM' \to \cB'$, which can be written down
explicitly (it involves Bessel functions, but nothing worse).  This is worked out in detail in \cite{Gaiotto:2008cd}.

Moreover, $g$ admits an explicit smooth extension
to a fibration $\cM \to \cB$, where $\cM \setminus \cM'$ consists of the fiber over $u=0$, a nodal torus.  This extended $g$
coincides with the well-known ``Ooguri-Vafa metric,'' first written down in \cite{Ooguri:1996me}.  So in this case our construction is a
new picture of the \hk structure on this known space.

One important drawback of this example is that it is only \ti{local} --- it is incomplete thanks to the boundary
of the disc $\cB$, and (as far as I know) has no suitable extension beyond this boundary.  This drawback is eliminated in
more interesting examples.
On the other hand this example is extremely simple and computable, thanks to the fact that the $\gamma$ for which
$\Omega(\gamma;u) \neq 0$ generate an isotropic lattice for $\IP{,}$.  Sadly, this virtue is also eliminated
in more interesting examples.

\section{More general singular loci}

In more interesting examples we cannot so easily study the behavior of the metric on $\cM'$
near the singular loci on $\cB$.
Nevertheless, we expect that the Ooguri-Vafa metric just discussed gives
a kind of local model for what happens generally near
the most generic kind of singular locus.  Namely, consider
some component $D_0 \subset D$, where
\begin{itemize}
\item $Z_{\gamma_0}(u) \to 0$ for some
specific $\gamma_0$,
\item $\Omega(\gamma_0;u) = 1$ for all $u$ in a neighborhood of
$D_0$,
\item $\gamma_0$ is primitive (i.e. there exists some $\gamma'$ with $\IP{\gamma_0, \gamma'} = 1$),
\item the monodromy of $\Gamma$ around $D_0$ is of ``Picard-Lefshetz type'', i.e.
\begin{equation}
 \gamma \to \gamma + \IP{\gamma, \gamma_0} \gamma_0.
\end{equation}
\end{itemize}
The most essential difference between this situation and the Ooguri-Vafa metric we
just discussed is that we no longer
require that $\Omega(\gamma;u) = 0$ for all $\gamma \neq \pm \gamma_0$.  Still,
near $D_0$ and for large enough $R$,
the quantum corrections coming from the charge $\gamma_0$, with $\Omega(\gamma_0;u) = 1$,
should dominate all others, and so $g$
should become similar to the Ooguri-Vafa metric.
In particular, at least for large enough $R$, $g$ should admit a smooth extension over $D_0$.
This remains to be rigorously understood.
I emphasize that it depends crucially on the condition $\Omega(\gamma_0;u) = 1$;
otherwise we would have no reason (either mathematical or physical) to expect
such a smooth extension of $g$ to exist.

All of the above admits an extension to the case where $\gamma_0$ is not primitive, but rather is
$k$ times a primitive vector.  In this case, instead of being smooth, we expect that
the completed $(\cM, g)$ has some mild singularities:  there should be
$k$ orbifold singularities of type $A_{k-1}$ lying over $D_0$.
This is still a significant improvement over the behavior of $g^\sf$.

The behavior of $g$ near higher-codimension strata on $D$ is more mysterious and should be very interesting.
At the moment it is not clear (at least to me) how to use our construction to get really
new information about it.

\section{Pentagon} \label{sec:pentagon}

The next simplest example is already much more nontrivial.
We fix a constant $\Lambda \in \C^\times$ (which enters the story in a trivial way:  it is safe to fix $\Lambda = 1$
if you prefer.)

\begin{vardescription}{Data 6}
\item[Data 1] $\cB$ is the complex plane, coordinatized by $u$.
\item[Data 2] The discriminant locus is $D = \{u = \pm 2 \Lambda^3 \} \subset \cB$.
So $\cB'$ is the twice-punctured plane.
\item[Data 3] Introduce a family of complex curves
\begin{equation}
\Sigma_u = \{ y^2 = z^3 - 3 \Lambda^2 z + u \} \subset \C^2.
\end{equation}
For $u \in \cB'$, $\Sigma_u$ is a noncompact smooth genus $1$ curve.
Define $\Gamma_u = H_1(\Sigma_u, \Z)$.
$\Gamma_u$ is a rank 2 lattice, the fiber of a local system $\Gamma$ over $\cB'$.  It is equipped with the intersection pairing $\IP{,}$.
$\Gamma_\flav$ is trivial.
\item[Data 4] Introduce the 1-form $\lambda = y\,\de z$.  $\lambda$ is a holomorphic 1-form on $\Sigma_u$, which would be meromorphic if extended
to the compactification of $\Sigma_u$ (it has a pole of order $6$ at the point at infinity, with zero residue).  Then for $\gamma \in \Gamma_u$,
\begin{equation}
 Z(\gamma) = \frac{1}{\pi} \oint_\gamma \lambda.
\end{equation}
\item[Data 5] Since $\Gamma_\flav$ is trivial, $\theta_\flav$ is trivial.
\item[Data 6] $\cB$ is divided into two domains $\cB_{\In}$ and $\cB_{\Out}$ (also sometimes called ``strong coupling'' and ``weak coupling'' respectively)
by the locus
\begin{equation}
W = \{ u: Z(\Gamma_u) \text{ is contained in a line in } \C \} \subset \cB.
\end{equation}
See Figure \ref{fig:pentagon-u-plane}.
\insfig{2.4in}{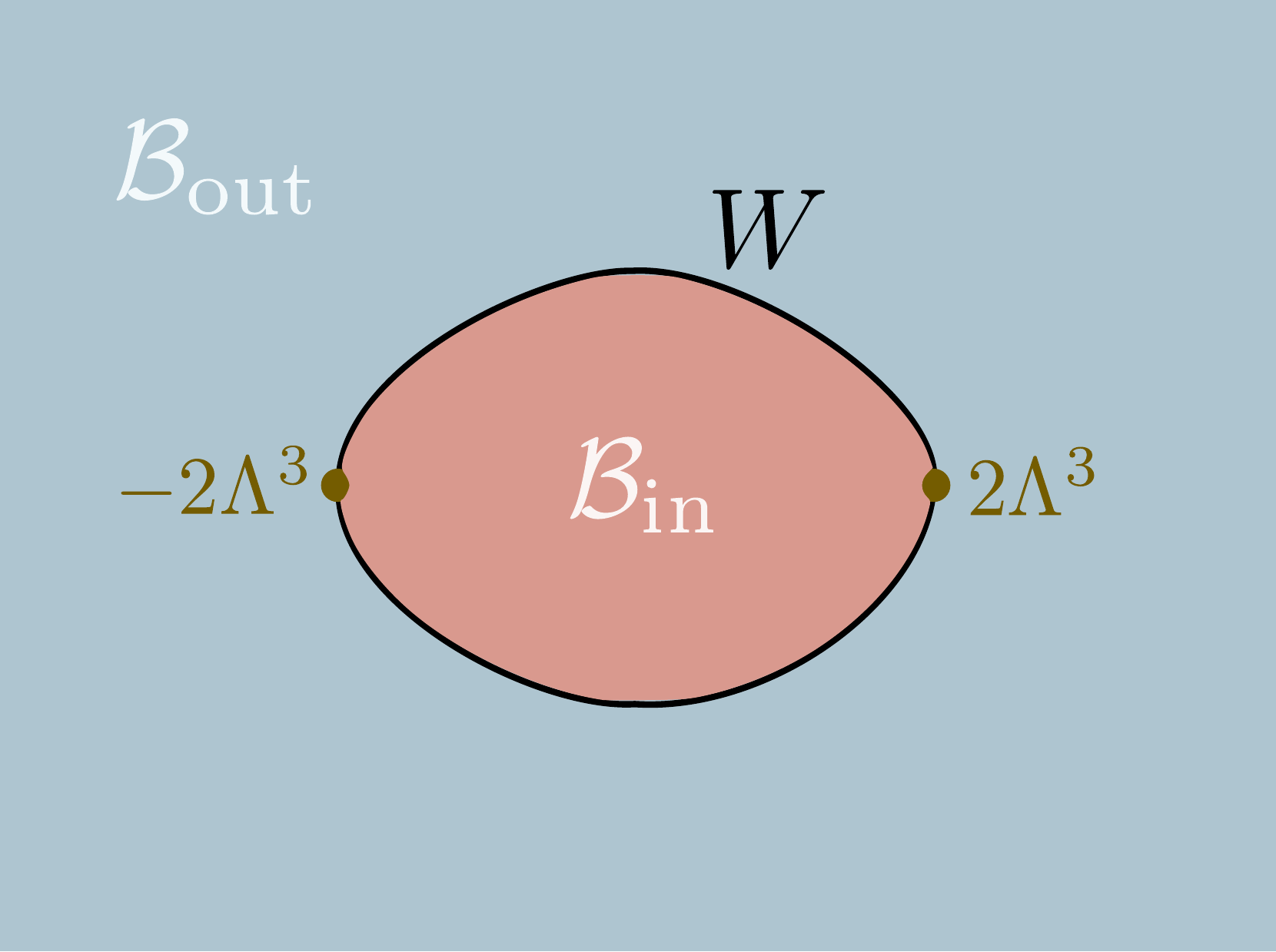}{The space $\cB$ in the example of Section \ref{sec:pentagon},
divided into two chambers by a wall.}
Since $\cB_\In$ is simply connected we may trivialize $\Gamma$ over $\cB_\In$ by primitive cycles
$\gamma_1$, $\gamma_2$ which collapse at the two points of $D$.  We choose them so that
$\IP{\gamma_1, \gamma_2} = 1$.
The set $\{\gamma_1, \gamma_2\}$ does not extend to a global trivialization of $\Gamma$, since it is not invariant under monodromy.
However, the set $\{ \gamma_1, \gamma_2, \gamma_1 + \gamma_2, -\gamma_1, -\gamma_2, -\gamma_1-\gamma_2\}$ is invariant under the monodromy
around infinity.  Therefore the following definition of $\Omega$ makes global sense:
\begin{align}
\text{For } u \in \cB_\In, \, \Omega(\gamma;u) &= \begin{cases} 1 \text{ for } \gamma \in \{\gamma_1, -\gamma_1, \gamma_2, -\gamma_2\}, \\ 0 \text{ otherwise.} \end{cases}\\
\text{For } u \in \cB_\Out, \, \Omega(\gamma;u) &= \begin{cases} 1 \text{ for } \gamma \in \{\gamma_1, -\gamma_1, \gamma_2, -\gamma_2, \gamma_1+\gamma_2, -\gamma_1-\gamma_2 \}, \\ 0 \text{ otherwise.} \end{cases}
\end{align}
\end{vardescription}
All of our conditions on the data are more or less trivial to check.  The most interesting one is the wall-crossing formula (Condition 7).
Here the question is:
choosing $u_{\In,\Out}$ to be two nearby points on opposite sides of $W$, and choosing $V$ to be a narrow sector which contains the rays $\ell_{\gamma_1}(u)$
and $\ell_{\gamma_2}(u)$ both for $u = u_\In$ and for $u = u_\Out$, do we have
\begin{equation}
A_V(u_\In) = \cK_{\gamma_1} \cK_{\gamma_2} \stackrel{?}{=} \cK_{\gamma_2} \cK_{\gamma_1 + \gamma_2} \cK_{\gamma_1} = A_V(u_\Out).
\end{equation}
This identity is indeed true:  it is the ``pentagon identity'' given in \cite{ks1}.  This identity can easily be checked by hand.  We remark in passing that (as also noted in \cite{ks1})
this identity is also closely related to the five-term
identity of the dilogarithm function and its quantum counterpart (see e.g. \cite{Faddeev:1993rs,MR2290758}.)

This example has the virtue that for every $u$ only finitely many $\Omega(\gamma;u)$ are nonvanishing.  This may lead to
some technical simplifications (although we emphasize that there should be no essential difference between this case and
the case where there are infinitely many nonvanishing $\Omega(\gamma;u)$, so long as the $\Omega(\gamma;u)$ grow slowly
enough with $\gamma$).

As we commented in the previous section, we expect that the metric $g$ on $\cM'$ in fact extends to a \ti{complete} metric on a space $\cM$,
obtained from $\cM'$ by adding nodal torus fibers over the two points of $D$, and that the metric around either of these nodal fibers looks like
the Ooguri-Vafa metric.

We believe that this complete metric actually has another name:  it is the metric on
a certain moduli space of rank-2 Higgs bundles on $\C\PP^1$
with an irregular singularity at $\infty$.  This point of view is discussed at some length in \cite{Gaiotto:2009hg}.
Also, for any $\zeta \in \C^\times$, the complex manifold $(\cM, J(\zeta))$ is isomorphic to
a partial compactification $\cM_{0,5}^{\mathrm{cyc}}$ of $\cM_{0,5}$, consisting of 5-tuples of points $(z_1, \dots, z_5)$ on $\C\PP^1$
where $z_i \neq z_{i+1}$ (with $i$ taken mod $5$).  Our description of this space is then closely related to
the discussion in \cite{qd-pentagon}.

\section{Hitchin systems} \label{sec:hitchin}

Finally I briefly describe a more geometric family of examples, considered in \cite{Gaiotto:2009hg}.

Fix a compact complex smooth curve $\bar{C}$.  Fix $n > 0$ marked points
$z_i \in \bar{C}$, and let $C = \bar{C} \setminus \{z_1, \dots, z_n\}$.
Also fix parameters $m_i \in \C$ and $m_i^{(3)} \in \R / 2 \pi \Z$ associated to the marked points.
Assume the $m_i$ and $m_i^{(3)}$ generic (in particular, the $m_i$ should be linearly independent over $\Q$.)

\begin{vardescription}{Data 6}
\item[Data 1] $\cB$ is the space of meromorphic quadratic differentials $\phi_2$ on $\bar{C}$ with double poles at each $z_i$,
of residue $m_i^2$.  (So $\cB$ is a complex affine space.)  To stay consistent with our previous notation we will use
either $u$ or $\phi_2$ to denote a point of $\cB$.
\item[Data 2] $D \subset \cB$ is the locus of $\phi_2$ which have at least one non-simple zero.  So $\cB'$ is the locus
of $\phi_2$ having only simple zeroes.
\item[Data 3] Let $T^*C$ be the holomorphic cotangent bundle to $C$.  For any fixed $u \in \cB$,
consider the noncompact complex curve
\begin{equation}
\Sigma_u = \{ (z \in C, \lambda \in T^*_z C):  \lambda^2 = \phi_2(z) \} \subset T^* C.
\end{equation}
For $u \in \cB'$, $\Sigma_u$ is smooth.  The obvious projection $\pi: \Sigma_u \to C$ is a double covering, branched over the
zeroes of $\phi_2$.
$\Sigma_u$ has a natural compactification $\bar \Sigma_u$ with a projection $\bar\pi: \bar\Sigma_u \to \bar C$.

$\Sigma_u$ is equipped with the involution $\lambda \mapsto -\lambda$.  Define $\Gamma_u$ to be the subgroup of $H_1(\Sigma_u, \Z)$ odd under this involution.
$\Gamma_u$ is the fiber of a local system $\Gamma$ over $\cB'$.  It is equipped with the intersection pairing $\IP{,}$.
$\Gamma_\flav \subset \Gamma$ is the radical of the pairing $\IP{,}$, which has rank $n$.
This radical does not undergo any monodromy as we vary $u$,
so we can think of $\Gamma_\flav$ as a single fixed lattice rather than a local system.
Finally, $\Gamma_\gauge = \Gamma / \Gamma_\flav$.

\item[Data 4]
By slight abuse of notation let $\lambda$ denote the Liouville (tautological) 1-form on $T^* C$.

Then for $\gamma \in \Gamma_u$, define
\begin{equation}
 Z(\gamma) = \frac{1}{\pi} \oint_\gamma \lambda.
\end{equation}

\item[Data 5]
The 1-form $\lambda$ restricted to $\Sigma_u$ extends meromorphically to $\bar\Sigma_u$, with
simple poles at the two preimages of $z_i$; let $z_i^\pm \in \bar\Sigma_u$ denote the preimage
at which $\lambda$ has residue $\pm m_i$.
The lattice $\Gamma_\flav$ has one generator $\gamma_{i,\flav}$ for each puncture $z_i$, given by the sum of a counterclockwise
loop around $z^+_i$ and a clockwise loop around $z^-_i$.
We define $\theta_\flav$ by
\begin{equation}
 \theta_\flav (\gamma_{i, \flav}) = m_i^{(3)}.
\end{equation}

\item[Data 6] The invariants $\Omega(\gamma;u)$ are defined in terms of the quadratic differential $\phi_2$,
as follows.

For any $\vartheta \in \R / 2 \pi \Z$,
define a \ti{$\vartheta$-trajectory} of $\phi_2$ to be a real curve $c \subset C$ such that, for any real tangent vector
$v$ to $c$, $\phi_2(v \otimes v) \in e^{2 \I \vartheta} \R_+$; call a $\vartheta$-trajectory \ti{maximal} if it is not properly contained
in any other $\vartheta$-trajectory.
The maximal $\vartheta$-trajectories make up a singular foliation of $C$, with three-pronged
singularities at the zeroes of $\phi_2$.

Define the \ti{mass} of a maximal $\vartheta$-trajectory $c$ to be $\int_c \abs{\sqrt{\phi_2}}$.  A generic maximal
$\vartheta$-trajectory has infinite mass; we are interested in the exceptional trajectories which have finite mass.
Let a \ti{finite $\vartheta$-trajectory} be a maximal $\vartheta$-trajectory with finite mass,
and a \ti{finite trajectory} be a pair $(c, \vartheta)$ where $c$ is a finite $\vartheta$-trajectory.
Finite trajectories come in two types:

\begin{itemize}
 \item \ti{Saddle connections}:  these are finite trajectories $c$ which ``run from one zero of $\phi_2$ to another,'' i.e.,
their boundary $\bar{c} \setminus c$ consists of two points (which are then necessarily zeroes of $\phi_2$).

 \item \ti{Closed loops}:  these are finite trajectories $c$ with the topology of $S^1$.  When such a trajectory occurs it sits in
a 1-parameter family of such trajectories, sweeping out an open annulus on $C$.
\end{itemize}

Given a finite trajectory $(c,\vartheta)$, define its \ti{lift} $\ell(c,\vartheta)$
to be the closure of $\pi^{-1}(c)$ on $\Sigma$.
$\ell(c,\vartheta)$ has no boundary; it is a single loop if $(c,\vartheta)$ is a saddle connection
(it is enlightening to draw a picture to see why),
and the disjoint union of two loops if $(c, \vartheta)$ is a closed loop.
The 1-form $e^{-i \vartheta} \lambda$ is real and nonvanishing on $\ell(c,\vartheta)$;
hence it induces an orientation on $\ell(c,\vartheta)$.
Note that if $(c,\vartheta)$ is a finite trajectory then $(c, \vartheta+\pi)$ is as well,
and $\ell(c,\vartheta)$ differs from $\ell(c,\vartheta+\pi)$ only by orientation reversal.
By construction, $\ell(c,\vartheta)$ is invariant under the combination of the deck
transformation $\lambda \mapsto -\lambda$ and orientation reversal.

For any $\gamma \in \Gamma_u$, let $SC(\gamma;u)$ be the set of all saddle connections $(c, \vartheta)$ with
$[\ell(c, \vartheta)] = \gamma$, and let $CL(\gamma;u)$ be the set of all isotopy classes of closed loops $(c, \vartheta)$
with $[\ell(c, \vartheta)] = \gamma$.
Now finally we can define
\begin{equation}
 \Omega(\gamma; u) = \# SC(\gamma;u) \, \, -2 \# CL(\gamma;u).
\end{equation}
(The strange-looking coefficients $+1$ and $-2$ here are really necessary --- otherwise the wall-crossing formula
(Condition 7) would not be satisfied!)

\end{vardescription}

These data satisfy all of our Conditions 1-7.  The most difficult to see are the last two.
Condition 6 follows from known results on quadratic differentials \cite{MR1053805,MR1827113} which say
$\Omega(\gamma;u)$ grows at most quadratically as a function of the coefficients of $\gamma$.
The wall-crossing formula (Condition 7) follows from a sort of inversion of the
logic we have followed up to this point:  namely, below we will give a direct description of
the complex spaces $(\cM, J(\zeta))$ and the functions $\cX_\gamma(x,\zeta)$ thereon
which solve the Riemann-Hilbert problem and are continuous except at the BPS rays.
The existence of such functions $\cX_\gamma(x,\zeta)$ then implies the wall-crossing formula
(following the discussion of Section \ref{sec:no-jumps}).

In \cite{Gaiotto:2009hg} we argued that the \hk space $\cM$ in this example is a
space of solutions of Hitchin equations on $\bar{C}$, with gauge group $PSU(2)$, and with ramification
at the marked points $z_i$ (with semisimple residues).  This is a much-studied space, considered in particular
in \cite{MR89a:32021,MR965220,MR887285,hbnc}.  In particular, it is known that the complex spaces $(\cM, J(\zeta))$
are moduli spaces of $PSL(2,\C)$ connections on $C$, with fixed eigenvalues of monodromy around $z_i$, given by
$\mu_\pm = \exp (\pm 2 \pi \I (\zeta^{-1} m_i - m_i^{(3)} - \zeta \bar{m}_i))$.

The $\cX_\gamma(x,\zeta)$ in this example are essentially functions considered earlier by Fock-Goncharov
in \cite{MR2233852}, themselves complexifications of the ``shear coordinates'' familiar in Teichm\"uller theory.
The main issue in identifying the Fock-Goncharov coordinates with our $\cX_\gamma(x,\zeta)$
is to prove that the Fock-Goncharov coordinates have the correct asymptotic behavior as $\zeta \to 0,\infty$.
This is accomplished by applying the WKB approximation
to a family of flat connections on $C$ of the form $\nabla(\zeta) = \varphi / \zeta + D + \bar\varphi \zeta$.

There is a generalization of this story to encompass quadratic differentials with
poles of order greater than $2$, also considered in  \cite{Gaiotto:2009hg}.
This generalization in particular includes the ``pentagon'' example
of Section \ref{sec:pentagon}; it corresponds to considering
quadratic differentials
$\varphi_2 = (z^3 - 3 \Lambda^2 z + u) \de z^2$ on $\C\PP^1$, with order-$7$ poles
at $z = \infty$.

Finally, we have extended many aspects of this story to the case of Hitchin equations
with higher rank gauge group $PSU(K)$ \cite{Gaiotto2012}.  In this case the coordinate functions $\cX_\gamma$ involve more
general coordinate systems than those which were described explicitly by Fock-Goncharov in \cite{MR2233852};
conjecturally the $\cX_\gamma$ exhaust the set of cluster coordinate systems.

\section{DT invariants}

Finally let us briefly consider another viewpoint on this story, which is really where it began.
The physical perspective on our construction makes clear that it should be
closely related to the theory of generalized
Donaldson-Thomas invariants (henceforth just ``DT invariants.'')  In this section
I briefly sketch that relation and a few examples.

\subsection{The dictionary}

In the theory of DT invariants, one begins with a triangulated category $\cD$
and constructs the space $\Stab(\cD)$ of \ti{Bridgeland stability conditions} on $\cD$ \cite{bridstab}.
Under some further conditions\footnote{which I am unfortunately not
competent to summarize} on $\cD$,
one then expects to be able to construct DT invariants depending on a point of $\Stab(\cD)$ \cite{Joyce:2008pc,ks1}, whose
dependence on the point of $\Stab(\cD)$ is governed by the wall-crossing formula.
In what follows I assume some familiarity with this story and formulate
the expected dictionary between the \hk data in our construction and the theory of DT invariants.
Many aspects of this dictionary are also described in Section 2.7 of \cite{ks1}.

\medskip

We need a few technical preliminaries to ``harmonize'' the two sides first:

\begin{itemize}
\item On the \hk data side:  suppose given an example of our Data 1-6 obeying our Conditions 1-7.
Fix a basepoint $u_0 \in \cB'$.  Let $\tilde \cB'$ denote the universal cover of $\cB'$.
Over this cover we may globally trivialize
the local system $\Gamma$, thus identifying all of its fibers with $\Gamma_{u_0}$.  The fiberwise homomorphism $Z: \Gamma \to \C$
can thus be thought of as a family of homomorphisms from the fixed lattice $\Gamma_{u_0}$ to $\C$,
depending on a point
$\tilde{u} \in \tilde{\cB'}$,
\begin{equation}
Z(\tilde{u}):  \Gamma_{u_0} \to \C.
\end{equation}

\item On the DT theory side:  suppose given an appropriate category $\cD$.
$\Stab(\cD)$ is a complex Poisson manifold,
carrying a natural ``forgetful'' map to $\Hom(K(\cD), \C)$
which is a local Poisson isomorphism  \cite{bridstab}.
We will consider a single connected component $\Stab^0(\cD) \subset \Stab(\cD)$.

\end{itemize}

We then have the following expected dictionary:

\medskip
\begin{center}
\begin{tabular}{|c|c|}
\hline
DT theory & \hk data \\
\hline
\hline
$K(\cD)$ & $\Gamma_{u_0}$ \\
Euler pairing & $\IP{\cdot,\cdot}$ \\
stability functions $Z: K(\cD) \to \C$ & $Z(\tilde{u}): \Gamma_{u_0} \to \C$ \\
DT invariants of $\cD$ & $c(\gamma) \in \Q$ from \eqref{eq:dt-inv} \\
a quotient of a Lagrangian $L \subset \Stab^0(\cD)$ & $\cB'$  \\
??? & $\cB$ \\
??? & $\theta_\flav$ \\
\hline
\end{tabular}
\end{center}

\bigskip

This dictionary has one especially awkward feature:  starting from the category $\cD$ it is not at all clear
how to choose the complex Lagrangian submanifold $L$.  Because of this problem, at the moment we do not
really have a recipe which begins with $\cD$ alone and constructs a corresponding
\hk space.  In particular examples which we do understand,
$L$ always has some nice geometric meaning (see the next section).
It would be very interesting to understand how to get $L$ in a purely categorical way.

\subsection{Examples}

For many examples of our construction of \hk metrics (probably in
all the examples that come from an underlying supersymmetric quantum
field theory, which includes all of the examples discussed so far
in this note), we expect that there is some triangulated
category $\cD$, fitting into the above dictionary.
Let us now describe a few examples:

\begin{itemize}

\item Let $\cD$ be the category of finite-dimensional modules over the Ginzburg dg algebra of the $A_2$ quiver.
In recent work of Sutherland \cite{Sutherland2011}, one connected component $\Stab^0(\cD) \subset \Stab(\cD)$ is identified with the universal cover of the
total space of a particular $\C^\times$ bundle over the moduli space $\cM_{1,1}$ of elliptic curves.
This result fits well into the above dictionary:  indeed we claim that the \hk data corresponding
to the category $\cD$ is that of the ``pentagon'' example of Section \ref{sec:pentagon} above.
The elliptic curves appearing in Sutherland's picture are the curves $\Sigma_u$
of Section \ref{sec:pentagon}.

\item Recent work of Bridgeland and Smith \cite{bridgelandsmith} is also relevant to this dictionary.

Begin with a real compact 2-manifold $C$, with $n>1$ marked points.
From the combinatorics of ideal triangulations of the curve $C$, one can build an associated
quiver $Q(C)$, using a superpotential function
first written down by Labardini-Fragoso \cite{Labardini-Fragoso2008}.\footnote{This quiver
and superpotential also appeared in the physics literature \cite{Alim:2011ae}.}
Let $\cD(C)$ be the derived category
of finite-dimensional modules over the Ginzburg dg algebra of $Q(C)$.
Bridgeland and Smith show
(roughly --- for the precise statement see \cite{bridgelandsmith}) that
there is a component $\Stab^0(\cD(C)) \subset \Stab(\cD(C))$, such that
a point of $\Stab^0(\cD(C))$ corresponds to a choice of complex structure on $C$
and a meromorphic quadratic differential thereon, with double poles at the marked points.
Among other things, this provides a family of nontrivial examples of categories $\cD$
where one has a geometric interpretation for at least a component of $\Stab(\cD)$.

This result fits in well with the dictionary proposed above:  it is consistent
with the idea that the categories $\cD(C)$ correspond to the
\hk data described in Section \ref{sec:hitchin}.
Moreover, revisiting Section \ref{sec:hitchin} we see that the
mysterious Lagrangian subspace $L \subset \Stab(\cD(C))$ appearing in the dictionary
has a nice meaning here:  it corresponds to fixing a particular
complex structure on $C$ and a choice of residues
at the marked points on $C$.

Bridgeland and Smith also consider a generalization corresponding
to allowing meromorphic quadratic differentials with higher-order poles.  This
generalization in particular gives another proof of Sutherland's results from \cite{Sutherland2011} (by considering quadratic differentials on $\C\PP^1$ with a single
pole of order $7$).

\item More ambitiously, at least on physical grounds
we expect that given a complete non-compact Calabi-Yau threefold $X$, both sides of this dictionary should exist.
Roughly speaking, $\cD = \cD(X)$ should be an appropriate version of the Fukaya category of $X$;
$\cB$ should be the moduli space of complex structures in $X$; $\Gamma$ should be
$H_3(X,\Z)$; $Z$ should be the period map; $c(\gamma)$ should be DT invariants counting
special Lagrangian 3-cycles in $X$.  The Lagrangian submanifold $L$ is the period domain of $X$;
the fact that it is Lagrangian is essentially Griffiths transversality.
Finally, the \hk space $\cM$ built by our construction is some version of the
family of intermediate Jacobians of $X$ (I say ``some version'' because we are dealing with non-compact
$X$).

The examples studied by Bridgeland and Smith, i.e. the examples of Section \ref{sec:hitchin} above,
also fall into this class.  The Calabi-Yau threefold $X(C)$ in this case
is a conic bundle over the curve $C$,
appropriately modified at the marked points; the
Fukaya category $\cD(X(C))$ is equivalent to the category $\cD(C)$ mentioned above.
This equivalence will also be explained in upcoming
work of Bridgeland and Smith.

(Incidentally, following through our various claims about the \hk space $\cM$ in these examples,
we see that on the one hand $\cM$ should be a version of the family of intermediate Jacobians of $X(C)$,
while on the other hand $\cM$ should be the $PSU(2)$ Hitchin system on $C$, with ramification at
the marked points.  So our claims would imply that these two integrable systems the same.
This equivalence is not really novel:  a version of it without the marked points
was described in \cite{Diaconescu:2006ry}, a fact which gives us some additional confidence in
our whole picture.)

For general $X$, it is not clear \ti{a priori} that the DT invariants will
grow slowly enough to satisfy our Condition 6.  Nevertheless, on physical grounds
we would expect the \hk manifold $\cM$ to exist for general $X$.\footnote{The idea is
that $\cM$ is the moduli space of the IIB string theory formulated on the 10-manifold
$X  \times S^1 \times \R^{2,1}$.}  Thus we expect that \ti{either} the DT invariants
do in fact grow slowly enough for us to prove that the Riemann-Hilbert problem
has a solution, \ti{or} they grow more quickly but
have some hidden extra structure that allows
the Riemann-Hilbert problem to have a solution anyway.

\item Finally let me describe a \ti{non}-example.
It is natural to ask:  what if we let $\cD$ be the Fukaya category of a
\ti{compact} Calabi-Yau threefold $X$ --- will there be corresponding \hk data then?
It seems that the answer is ``yes'' --- we can define the data by the same recipe as we use
for non-compact $X$ --- but these data would
not satisfy precisely our Conditions 1-7.
In particular, Condition 4 (positive definiteness) will certainly be violated.  However, this
violation is of a rather controlled sort;
there is just one negative direction.  So, were this the only difficulty, the
expected consequence would be that the space $\cM$ we obtain is not \hk but pseudo-\hk, with one
negative direction.  ($\cM$ in this case is the family of intermediate
Jacobians of $X$, fibered over the moduli space of polarized complex structures on $X$.
These intermediate Jacobians are quotients of $H^{3,0} \oplus H^{2,1} (X)$,
and the negative direction is coming from $H^{3,0}(X)$;
it is related to the fact that when
equipped with its ``Griffiths'' complex structure, the intermediate Jacobian
is not principally polarized.)  However, there is also a second, more serious difficulty:  the
invariants $\Omega(\gamma)$ counting special Lagrangian 3-cycles in $X$ are expected
to grow very quickly as functions of $\gamma$ (roughly $\Omega(\gamma) \sim \exp c \norm{\gamma}^2$),
badly violating our Condition 6.  As a result
it is far from clear whether our construction of \hk metrics
should be directly applicable to this situation.

This difficulty is in some sense anticipated in the physics literature.  Indeed,
physics does \ti{not} predict directly that there is an \hk manifold associated to a compact Calabi-Yau
threefold $X$.
Rather it predicts the existence of a \ti{\qk} manifold.
As in the \hk case, it should be possible to construct the desired \qk structure by beginning
with a simple ``semi-flat'' metric $g^\sf$ and modifying
it by quantum corrections.\footnote{In this case ``semi-flat'' means that $g^\sf$ is locally
invariant under a Heisenberg group of isometries, replacing the torus group that appeared in the \hk
case.}
The semi-flat metric in this case was first described by
Ferrara and Sabharwal in
\cite{Ferrara:1989ik}, and was recently discussed by Hitchin in \cite{MR2494168}.
The description of the quantum corrections has been studied intensely in physics, with various
interesting partial results.  In particular, some of the quantum
corrections are expected to be precise analogues of the ones
we have described in the \hk case, indeed related by a ``\qk/\hk correspondence''
\cite{Alexandrov:2011ac,MR2394039}.  However, one also expects new quantum
corrections in the \qk case which do not have an \hk analogue.  As far as I know,
there are no examples yet of $X$ where all quantum corrections have been fully described.

\end{itemize}

\bibliography{hyperkahler-notes}

\small\normalsize

\end{document}